\documentclass[11pt,a4j]{article}
\usepackage{amsthm}
    \newtheorem*{inthm}{Theorem}
    
    \newtheorem{thm}{Theorem}[section]
    \newtheorem{prop}[thm]{Proposition}
    \newtheorem{lem}[thm]{Lemma}
    \newtheorem{cor}[thm]{Corollary}
    
  \theoremstyle{definition}
    \newtheorem*{que}{Question}
    \newtheorem{defi}[thm]{Definition}
  \theoremstyle{remark}
    \newtheorem{rem}[thm]{Remark}
    \newtheorem{ex}[thm]{Example}

\usepackage{amsmath}
\usepackage{amssymb}
\usepackage[all]{xy}
\usepackage[dvipdfmx]{graphicx}
\usepackage{mathrsfs}
\usepackage{bm}

\title{Developing Maps and Engel Diffeomorphisms}
\author{Koji Yamazaki}
\begin{document}
\maketitle	
\begin{abstract}
A completely nonintegrable $2$-dimensional distribution on a $4$-manifold is called an Engel structure.
A $4$-manifold with an Engel structure is called an Engel manifold.
The developing map for an Engel manifold is very important tool to determine the Engel structure.
Montgomery \cite{montgomery1993generic} used it to prove that an Engel automorphism is determined by the values on a global slice.
Moreover, Montgomery \cite{montgomery1993generic} constructed Engel manifolds whose automorphism group is small.
In this paper, we prove that the automorphism group of an Engel manifold is embedded into the automorphism group of the Cartan prolongation of a contact $3$-orbifold, if the developing map is not a covering map.
As an application, we will construct an Engel manifold whose automorphism group is trivial.
\end{abstract}
\setcounter{section}{-1}								
	\section{Introduction} \label{sec0}
An Engel structure is a completely nonintegrable $2$-dimensional distribution $\mathcal{D}$ on a $4$-manifold $E$, and the pair $(E, \mathcal{D})$ is called an Engel manifold
(cf. Definition \ref{defEngel}).
\par
In 1889, Engel \cite{Engel} proved a Darboux-type theorem for Engel manifolds, which states that every Engel manifold is locally Engel diffeomorphic to $\mathbb{R}^4$ with standard Engel structure.
In 1993, Montgomery \cite{montgomery1993generic} proved that every stable germ of $k$-dimensional distribution on a $n$-manifold must satisfy the inequality $k(n-k) \le n$.
(See \cite{montgomery1993generic}.)
A pair $(n, k)$ of positive integers satisfies $k(n-k) \le n$ and $0 < k < n$ if and only if it satisfies $k = 1$, $k = n-1$, or $(n, k) = (4, 2)$.
The equality $k = 1$ is realized by $1$-dimensional foliations, while the equality $k = n-1$ is realized by contact structures if $n$ is odd and even contact structures if $n$ is even.
The equality $(n, k) = (4, 2)$ is realized by Engel structures.\par
An {\it Engel automorphism} (resp. a {\it contact automorphism}) of an Engel manifold (resp. a contact manifold) is a self-diffeomorphism preserving the Engel structure (resp. the contact structure).
The automorphism group $\operatorname{Aut}(E, \mathcal{D})$ (resp. $\operatorname{Aut}(M, \xi)$) is defined as the set of Engel automorphisms (resp. contact automorphisms) of an Engel manifold $(E, \mathcal{D})$ (resp. a contact manifold $(M, \xi)$).
An Engel manifold has a $1$-dimensional distribution called the {\it characteristic foliation} (cf. Definition \ref{defchar}).
Let $(E, \mathcal{D})$ be an Engel manifold, and let $M \subset E$ be a $3$-dimensional submanifold intersecting transversally with the characteristic foliation of the Engel manifold $(E, \mathcal{D})$.
(Such a submanifold is called a {\it slice} of a foliation.)
Then, the manifold $M$ has a contact structure $TM \cap \mathcal{D}^2$.
(See \cite{yamazaki2018engel}.)
It is known that an Engel structure is completely determined locally by the characteristic foliation and the contact structure of a slice of the characteristic foliation.
\vspace{5pt}\par
Montgomery proved the following theorem.
\begin{inthm}[\mbox{Montgomery \cite[Theorem 3 (b)]{montgomery1999engel}}]
Let $(E, \mathcal{D})$ be an Engel manifold, and let $\Sigma \subset E$ be a global slice of the characteristic foliation $\mathcal{L}$. 
(i.e. $\Sigma$ is a $3$-dimensional submanifold intersecting transversally with all leaves of the foliation $\mathcal{L}$.)
Two Engel automorphisms $f$ and $g$ of $(E, \mathcal{D})$ satisfy an equality $f = g$ if and only if they satisfy an equality $f(x) = g(x)$ for any $x \in \Sigma$.
\end{inthm}
Montgomery's theorem may be useful for computing Engel automorphism groups.
In fact, Montgomery \cite[Theorem 6]{montgomery1999engel} used this theorem to discover an Engel manifold with small automorphism group.
Moreover, in practice, the theorem can be improved.
If we could construct a global slice functorially, we might obtain a homomorphism from an Engel automorphism group to a contact automorphism group.
Then, Montgomery's theorem should be rephrased as a claim that the group homomomorphism is injective.
However, there is no hope to construct of global slices functorially.
To modify this, we use the leaf space instead of a global slice. \par
Remark, however, that the leaf space is not a manifold in general.
We say that an Engel manifold has the {\it trivial characteristic foliation} if the leaf space of the characteristic foliation is a manifold.
\vspace{5pt}\par
Any (global) slice of the characteristic foliation $\mathcal{L}$ of an Engel manifold has a contact structure, and any holonomy of $\mathcal{L}$ preserve the contact structure.
This means the leaf space of $\mathcal{L}$ has a contact structure.
(See \cite{yamazaki2018engel}.)
We write $\mathbb{L}(E, \mathcal{D})$ for the pair of the leaf space of the characteristic foliation of an Engel manifold $(E, \mathcal{D})$ and the contact structure.
The leaf space $\mathbb{L}(E, \mathcal{D})$ is a contact $3$-dimesional manifold if $(E, \mathcal{D})$ has the trivial characteristic foliation
(cf. Proposition \ref{proplspman}).
The correspondence $(E, \mathcal{D}) \mapsto \mathbb{L}(E, \mathcal{D})$ induces a functor $\mathbb{L}$ from the category of Engel manifolds with trivial characteristic foliations to the category of contact $3$-dimensional manifolds
(cf. Section \ref{sec1.4}). \par
There is also an inverse construction. 
For a contact $3$-manifold $(M, \xi)$, the projectization $\mathbb{P}(\xi)$ of the $2$-plane bundle $\xi$ admits an Engel structure
(cf. Section \ref{sec1.1.3}).
This Engel manifold is called the {\it Cartan prolongation} of $(M, \xi)$, and is denoted by $\mathbb{P}(M, \xi)$.
The correspondence $(M, \xi) \mapsto \mathbb{P}(M, \xi)$ induces a functor $\mathbb{P}$ from the category of contact $3$-manifolds to the category of Engel manifolds with trivial characteristic foliations
(cf. Section \ref{sec1.4}).
The functor $\mathbb{P}$ is fully-faithful
(cf. Corollary \ref{corfullyfaithful} and Section \ref{sec1.4}). \par
Cartan prolongations can also be defined for a larger category, \'{e}tale Lie groupoids.
In this paper, we only need to discuss Cartan prolongations of contact orbifolds
(cf. Remark \ref{remnoneed}).
Although the leaf space $\mathbb{L}(E, \mathcal{D})$ of the characteristic foliation of an Engel manifold $(E, \mathcal{D})$ is not a manifold in general, the Cartan prolongation $\mathbb{P}\mathbb{L}(E, \mathcal{D})$ of $\mathbb{L}(E, \mathcal{D})$ can still be a manifold.
We say that an Engel manifold $(E, \mathcal{D})$ has the {\it good characteristic foliation} if the Cartan prolongation $\mathbb{P}\mathbb{L}(E, \mathcal{D})$ of the leaf space $\mathbb{L}(E, \mathcal{D})$ is a manifold
(cf. Definition \ref{goodchar}). \par
The composition of two functors $(E, \mathcal{D}) \mapsto \mathbb{L}(E, \mathcal{D})$ and $(M, \xi) \mapsto \mathbb{P}(M, \xi)$ induces a group homomorphism $\Phi = \mathbb{P} \circ \mathbb{L} : \operatorname{Aut}(E, \mathcal{D}) \rightarrow \operatorname{Aut}(\mathbb{L}(E, \mathcal{D})) \cong \operatorname{Aut}(\mathbb{P}\mathbb{L}(E, \mathcal{D})$ between the automorphism groups.
We want to find out conditions under which the group homomorphism $\Phi$ is injective.
\vspace{5pt}\par
Let $(E, \mathcal{D})$ be an Engel manifold.
(One can, suppose that the Engel manifold $(E, \mathcal{D})$ has the good characteristic foliation for simplicity.
See \cite{yamazaki2018engel} for the general case.)
There exists a local Engel diffeomorphism $\phi : (E, \mathcal{D}) \rightarrow \mathbb{P}\mathbb{L}(E, \mathcal{D})$ called the {\it developing map} associated with the Engel manifold $(E, \mathcal{D})$
(cf. Definition \ref{defdev}), where a {\it local Engel diffeomorphism} is a local diffeomorphism between Engel manifolds preserving the Engel structures
(cf. Definition \ref{defEngelmor}).
The developing map is natural.
i.e. For any Engel automorphism $f : (E, \mathcal{D}) \rightarrow (E, \mathcal{D})$, the following diagram is commutative:
	\[\xymatrix{
		(E, \mathcal{D}) \ar[r]^-{\phi}
		& \mathbb{P}\mathbb{L}(E, \mathcal{D})
	\\
		(E, \mathcal{D}) \ar[r]^-{\phi} \ar[u]^-{f}
		& \mathbb{P}\mathbb{L}(E, \mathcal{D}) \ar[u]^-{\mathbb{P}\mathbb{L}f}.
	}\]
(Where $\mathbb{P}\mathbb{L}f$ is an Engel automorphism $\Phi(f)$.
Recall that $\Phi$ is the group homomorphism $\mathbb{P} \circ \mathbb{L} : \operatorname{Aut}(E, \mathcal{D}) \rightarrow \operatorname{Aut}(\mathbb{P}\mathbb{L}(E, \mathcal{D})$.)\par
We now review an outline of the proof of Montgomery's theorem.
Let $\Sigma \subset E$ be a global slice of the characteristic foliation $\mathcal{L}$ of the Engel manifold $(E, \mathcal{D})$.
Let $f, g : (E, \mathcal{D}) \rightarrow (E, \mathcal{D})$ be two Engel automorphisms.
We suppose that $\mathcal{L}$ is trivial for simplicity.
The condition $f(x) = g(x)$ for any $x \in \Sigma$ implies an equality $\mathbb{L}f = \mathbb{L}g$.
Then, we obtain an equality $\mathbb{P}\mathbb{L}f = \mathbb{P}\mathbb{L}g$.
Let $\sigma$ be the Engel diffeomorphism $\mathbb{P}\mathbb{L}f( = \mathbb{P}\mathbb{L}g)$.
Fix a point $x \in \Sigma$.
Let $L$ be a leaf of $\mathcal{L}$ including the point $x$.
Let $L'$ be a leaf of the characteristic foliation of the Engel manifold $\mathbb{P}\mathbb{L}(E, \mathcal{D})$ including the point $\phi(x)$.
We discuss the following diagram:
	\[\xymatrix{
		\phi^{-1} (\operatorname{Im}(\sigma |_{L'})) \ar[r]^-{\phi}
		& \operatorname{Im}(\sigma |_{L'})
	\\
		L \ar@<0.5ex>[u]^-{f} \ar@<-0.5ex>[u]_-{g} \ar[r]_-{\phi}
		& L' \ar[u]^-{\sigma}.
	}\]
If it satisfies an equality $f(x) = g(x)$, then it satisfies an equality $f(y) = g(y)$ for any $y \in L$, because $L$ is connected and $f$ and $g$ are lifts of a continous map $\sigma \circ \phi : L \rightarrow \operatorname{Im}(\sigma |_{L'})$ with respect to the local diffeomorphism $\phi : \phi^{-1} (\operatorname{Im}(\sigma |_{L'})) \rightarrow \operatorname{Im}(\sigma |_{L'})$.
If a local diffeomorphism $\phi : L \rightarrow L'$ is not a covering map, then the condition $f(x) = g(x)$ is not necessary.
This is because, roughly speaking, a ``point" $\phi(\infty)$ is included in $L'$ in this case, where $\infty$ is an ``infinity point" of $L$
(cf. Lemma \ref{noncov}).
Then, it satisfies the equality $f(\infty) = g(\infty)$.\par
The main theorem in this paper is the following result.
\setcounter{section}{2}
\setcounter{thm}{0}
\begin{thm} \label{main}
Let $(E, \mathcal{D})$ be a connected Engel manifold with good characteristic foliation
(cf. Definition \ref{goodchar}).
If the developing map $\phi$ associated with $(E, \mathcal{D})$ is not a covering map, then the group homomorphism $\Phi = \mathbb{P} \circ \mathbb{L} : \operatorname{Aut}(E, \mathcal{D}) \rightarrow \operatorname{Aut}(\mathbb{P}\mathbb{L}(E, \mathcal{D}))$ is injective.
\end{thm}
In Montgomery's theorem, Engel automorphisms are determined by the values on a global slice. 
In Theorem \ref{main}, Engel automorphisms are determined by the values on the infinity points of the leaves.
Roughly speaking, Theorem \ref{main} is regarded as an application of Montgomery's Theorem to 
\[
\Sigma = \{ \mbox{Bounded infinity points of a leaf of the characteristic foliation} \}.
\]
An infinity point $\infty$ being {\it bounded} means that the point $\phi(\infty)$ is included in the Cartan prolongation $\mathbb{P}\mathbb{L}(E, \mathcal{D})$.
The condition that $\phi$ is not covering implies $\Sigma \neq \emptyset$.
\vspace{5pt}\par
Finally, as an application, we give a positive answer (Theorem \ref{app}) to the following question.
\begin{que}[AIM Problem lists\footnote{http://aimpl.org/engelstr/3/}]
Is there an Engel manifold with trivial automorphism group?
\end{que}
An example of Engel manifolds with trivial automorphism group is constructed in Example \ref{ex3}.
Our construction is very pathological.
However, our intuition predicts that the automorphism group of a generic Engel manifold is trivial, or very small.
Mitsumatsu gives the following question, which remains open.
\begin{que}[Mitsumatsu]
Is there a {\it closed} Engel manifold with trivial automorphism group?
\end{que}
In Section 1.1.1, we will give a definition of Engel manifolds.
In Section 1.1.2, we will define the characteristic foliation and discuss the functor from Engel to contact.
In Section 1.1.3, we will define the Cartan prolongation and discuss the functor from contact to Engel.
In Section 1.1.4, we will define the developing map.
In Section 1.2, we generalize the discussion in the previous sections.
We prove the main theorem in Section 2.2.
We give a positive answer to the AIM Problem in Section 2.3.
		\subsection*{Acknowledgements}
I would like to thank Professor Mitsumatsu and the members of the Saturday Seminar for their discussions.
\setcounter{section}{0}
\setcounter{thm}{0}
	\section{Engel manifolds} \label{sec1}
In Section 1.1, we review some basic concepts and their properties of Engel manifolds.
In Section 1.2, we introduce contact orbifolds to extend the relationship between Engel structures and contact structures.
In Section 1.3, we will reformulate the previous discussion in terms of category theory.
	\subsection{Engel manifolds} \label{sec1.1}
An Engel structure is a completely nonintegrable $2$-dimensional distribution $\mathcal{D}$ on a $4$-manifold $E$.
The pair $(E, \mathcal{D})$ is called an Engel manifold.
Any Engel manifold has a $1$-dimensional foliation called the {\it characteristic foliation}.
The leaf space of the characteristic foliation admits a canonical contact structure.
Conversely, a contact $3$-manifold has an Engel manifold called the {\it Cartan prolongation}.
An Engel structure determines the {\it developing map}. \par
This section is based on \cite{montgomery1999engel}, \cite{adachi2002engel} and \cite{vogel2009existence}.
The detailed proofs are summarized in \cite{yamazaki2018engel}.
	\subsubsection{Definitions of Engel manifolds}
\begin{defi}[Engel manifolds] \label{defEngel}
Let $E$ be a $4$-dimensional manifold.
An {\it Engel structure} on $E$ is a $2$-dimensional distribution $\mathcal{D} \subset TE$ such that subsets $\mathcal{D}_x^2 \subset T_xE$ indexed by points $x \in E$ determines $3$-dimensional distribution $\mathcal{D}^2$ and that subsets $\mathcal{D}_x^3 \subset T_xE$ indexed by points $x \in E$ determines $4$-dimensional distribution $\mathcal{D}^3$, where the sets $\mathcal{D}_x^2$ and $\mathcal{D}_x^3$ are defined as:
\[
\begin{array}{lcl}
\mathcal{D}_x^2 & = & \{ [X, Y]_x + Z_x \, | \, X_x, Y_x, Z_x \in \mathcal{D}_x \}, \\
\mathcal{D}_x^3 & = & \{ [X, Y]_x + Z_x \, | \, X_x, Y_x, Z_x \in \mathcal{D}_x^2 \}.
\end{array}
\]
The pair $(E, \mathcal{D})$ is called an {\it Engel manifold}.
\end{defi}
\begin{rem}
An Engel manifold can also be defined as follows.
Let $\mathcal{D}$ be a smooth distribution on a manifold.
We can regard $\mathcal{D}$ as a locally free sheaf of some vector fields.
Define the sheaf $\mathcal{D}^n$ as the sheafification of a presheaf $\mathcal{D'}^n$ where the presheaf $\mathcal{D'}^n$ is inductively defined as follows.
\[
\begin{array}{lcl}
\mathcal{D'}^1 & = &\mathcal{D} \\
\mathcal{D'}^{n+1} & = & [\mathcal{D}^n, \mathcal{D}^n] + \mathcal{D}^n \\
& = & \{ [X, Y] + Z \, | \, X, Y, Z \in \mathcal{D}^n \}.
\end{array}
\]
Suppose that $\mathcal{D}$ is a $2$-dimensional distribution on a $4$-manifold.
Then, $\mathcal{D}$ is an Engel structure if and only if the sheaf $\mathcal{D}^2$ is a locally free sheaf whose rank is $3$ and the sheaf $\mathcal{D}^3$ is a locally free sheaf whose rank is $4$.
\end{rem}
A morphism between Engel manifolds is defined as follows.
\begin{defi} \label{defEngelmor}
A {\it local diffeomorphism} $f : E_1 \rightarrow E_2$ is a smooth map $f : E_1 \rightarrow E_2$ such that the differential $df_x : T_xE_1 \rightarrow T_{f(x)}E_2$ is an isomorphism for any $x \in E_1$. \par
Let $(E_1, \mathcal{D}_1)$ and $(E_2, \mathcal{D}_2)$ be Engel manifolds. 
A {\it local Engel diffeomorphism} $f : (E_1, \mathcal{D}_1) \rightarrow (E_2, \mathcal{D}_2)$ is a local diffeomorphism $f : E_1 \rightarrow E_2$ with $df(\mathcal{D}_1) \subset \mathcal{D}_2$. \par
A local Engel diffeomorphism $f : (E_1, \mathcal{D}_1) \rightarrow (E_2, \mathcal{D}_2)$ is an {\it Engel diffeomorphism} if the map $f : E_1 \rightarrow E_2$ is a bijection.
An {\it Engel automorphism} is an Engel diffeomorphism from an Engel manifold to itself.
\end{defi}
For example, {\it developing maps} are local Engel diffeomorphisms
(cf. Section \ref{sec1.1.4}).
A morphism of the category of Engel manifolds is defined as a local Engel diffeomorphism
(cf. Section \ref{sec1.4}).
	\subsubsection{Characteristic foliations}
We will discuss a correspondence from Engel to contact.
The main subject of this subsubsection is the {\it characteristic foliation}.
Any Engel manifold has a $1$-dimensional distribution called the characteristic foliation.
(Remark that any $1$-dimensional distribution is integrable.
An integrable distribution determines a foliation.
In this paper, we do not distinguish an integrable distribution with the corresponding foliation unless it is confusing.)
\begin{prop}[\mbox{\cite{montgomery1999engel}, cf. \cite[section 1.2.1]{yamazaki2018engel}}]
Let $(E, \mathcal{D})$ be an Engel manifold.
Then, there exists a unique $1$-dimensional distribution $\mathcal{L}$ on the manifold $E$ satisfying $[\mathcal{D}^2, \mathcal{L}] \subset \mathcal{D}^2$.
Moreover, the distribution $\mathcal{L}$ included in $\mathcal{D}$.
\end{prop}
\begin{defi} \label{defchar}
The above $1$-dimensional distribution $\mathcal{L}$ is called the {\it characteristic foliation} of the Engel manifold $(E, \mathcal{D})$.
\end{defi}
The leaf space of the characteristic foliation has a contact structure.
(See \cite{yamazaki2018engel}.)
However, the leaf space may not be a manifold.
In order to describe the correspondence between Engel manifolds and contact manifolds, we give a definition of contact manifolds and trivial characteristic foliations.
\begin{defi}[contact manifolds]
Let $M$ be a $(2n+1)$-dimensional manifold.
A {\it contact structure} on the manifold $M$ is a corank $1$ distribution $\xi \subset TM$ such that a $(2n+1)$-form $(\wedge^n d \alpha) \wedge \alpha$ does not vanish anywhere for any local $1$-form $\alpha$ with $\xi = \operatorname{Ker}(\alpha)$.
The pair $(M, \xi)$ is called a {\it contact manifold}, and the above $1$-form $\alpha$ is called a {\it contact form}. \par
Let $(M_1, \xi_1)$ and $(M_2, \xi_2)$ be contact manifolds. 
A {\it local contactomorphism} $f : (M_1, \xi_1) \rightarrow (M_2, \xi_2)$ is a local diffeomorphism $f : M_1 \rightarrow M_2$ with $df(\xi_1) \subset \xi_2$.
A local contactomorphism $f : (M_1, \xi_1) \rightarrow (M_2, \xi_2)$ is a {\it contactomorphism} if the map $f : M_1 \rightarrow M_2$ is a bijection.
A {\it contact automorphism} is a contactomorphism from a contact manifold to itself.
\end{defi}
\begin{rem} \label{remlocal}
Any example of local contactomorphism does not appear in this paper.
However, it is a necessary definition to neatly describe categorical arguments (cf. Section \ref{sec1.4}), such as functor.
\end{rem}
\begin{defi} \label{trichar}
Let $(E, \mathcal{D})$ be an Engel manifold.
We say that $(E, \mathcal{D})$ has the {\it trivial characteristic foliation} if the leaf space of the characteristic foliation is a manifold.
\end{defi}
\begin{prop}[\mbox{\cite{montgomery1999engel}, cf. \cite[Proposition 1.15]{yamazaki2018engel}}] \label{proplspman}
Let $(E, \mathcal{D})$ be an Engel manifold with trivial characteristic foliation $\mathcal{L}$.
Let $M$ be the leaf space of the foiation $\mathcal{L}$.
Then, the distribution $\xi = \mathcal{D}^2 / \mathcal{L}$ is well-defined on $M$, and it is a contact structure on $M$.
\end{prop}
The pair $(M, \xi)$ of the above leaf space $M$ of $\mathcal{L}$ and the above contact structure $\xi$ is a contact manifold.
The contact manifold $(M, \xi)$ is denoted by $\mathbb{L}(E, \mathcal{D})$.
In this paper, we do not distinguish the contact manifold $\mathbb{L}(E, \mathcal{D})$ with the manifold $M$ unless it is confusing.
	\subsubsection{Cartan prolongations} \label{sec1.1.3}
We will construct an Engel manifold, called the {\it Cartan prolongation}, from a contact $3$-manifold.
Let $(M, \xi)$ be a contact 3-manifold.
Define a manifold $E$ as $E = \mathbb{P} (\xi) = \displaystyle{\coprod_{x \in M} } \mathbb{P} (\xi_x) \times \{ x \}$, where $\mathbb{P}(V)$ is a projectization of a vector space $V$.
(i.e. $\mathbb{P}(V) = \{ \mbox{dimension $1$ subspaces of $V$} \}$)
Let $\pi : E \rightarrow M$ be the projection. \par
We define a $2$-dimensional distribution $\mathcal{D}$ on $E$ in the following way. 
For each $(l, x) \in E$ with $\pi (l, x) = x$, $l \subset \xi_x$ is a $1$-dimensional linear subspace of $\xi_x$.
Define a distribution $\mathcal{D}$ as $\mathcal{D}_{(l, x)} = d \pi_{(l, x)}^{-1} (l) \subset T_{(l, x)} E$.
\begin{lem}[\mbox{\cite{montgomery1999engel}, cf. \cite[Lemma 1.16]{yamazaki2018engel}}] \label{Cartanlem}
\mbox{}
\begin{itemize}
\item The above distribution $\mathcal{D}$ is an Engel structure on $E$.
\item The distribution $\mathcal{D}^2$ coincides with the distribution $d\pi^{-1}(\xi)$.
\item The distribution $\operatorname{Ker}(d\pi)$ is the characteristic foliation of the Engel manifold $(E, \mathcal{D})$.
\end{itemize}
\end{lem}
\begin{defi}[Cartan prolongation]
The above Engel manifold $(E, \mathcal{D})$ is called the {\it Cartan prolongation} of the contact manifold $(M, \xi)$.
This Engel manifold $(E, \mathcal{D})$ is denoted by $\mathbb{P}(M, \xi)$.
In this paper, we do not distinguish the Engel manifold $\mathbb{P}(M, \xi)$ with the manifold $E$ unless it is confusing.
\end{defi}
	\subsubsection{Developing maps} \label{sec1.1.4}
The leaf space of the Cartan prolongation $\mathbb{P}(M, \xi)$ is contactomorphic to a given contact manifold $(M, \xi)$.
The Cartan prolongation is ``minimal" object among such Engel manifolds.
The ``minimality" means a certain universality.\par
Let $(E, \mathcal{D})$ be an Engel manifold with trivial characteristic foliation $\mathcal{L}$.
Let $(M, \xi) = \mathbb{L}(E, \mathcal{D})$, and let $\pi : E \rightarrow M$ be the projection.
Define $\phi : E \rightarrow \mathbb{P}(\xi)$ as $\phi(e) = d\phi(\mathcal{D}_e)$ for $e \in E$.
\begin{lem}[\mbox{\cite{montgomery1999engel}, cf. \cite[Lemma 1.18]{yamazaki2018engel}}]
The above map $\phi$ is a local Engel diffeomorphism.
\end{lem}
\begin{defi} \label{defdev}
The above local Engel diffeomorphism $\phi$ is called the {\it developing map} associated with the Engel manifold $(E, \mathcal{D})$.
\end{defi}
\begin{rem} \label{rembij1}
Let $(E, \mathcal{D})$ be an Engel manifold with trivial characteristic foliation $\mathcal{L}$.
Let $\mathcal{L'}$ be the characteristic foliation of the Engel manifold $\mathbb{P}\mathbb{L}(E, \mathcal{D})$.
The leaf space of the foliation $\mathcal{L'}$ is contactomorphic to the contact manifold $\mathbb{L}(E, \mathcal{D})$.
Then, the developing map $\phi$ induces the following bijection $\phi_\ast$.
\[
\{ \mbox{the leaves of the foliation $\mathcal{L}$} \} \overset{\phi_\ast}{\rightarrow} \{ \mbox{the leaves of the foliation $\mathcal{L'}$} \}.
\]
This holds even if the leaf space of the foliation $\mathcal{L}$ is not a manifold
(cf. Remark \ref{rembij2}).
\end{rem}
The developing map has a universal property.
(Remark that the correspondence $(M, \xi) \mapsto \mathbb{P}(M, \xi)$ induces C from B.)
\begin{prop}[\mbox{cf. \cite[Proposition 1.19]{yamazaki2018engel}}] \label{propuniv}
Take any contact $3$-manifold $(N , \xi_N)$ and any local Engel diffeomorphism $\psi : E \rightarrow \mathbb{P} (\xi_N)$.
Then, there exists a unique local contactomorphism $\tilde{\psi} : M \rightarrow N$ such that the following diagram is commutative:
	\[\xymatrix{
		E \ar[r]^-{\phi} \ar[rd]_-{\psi}
		& \mathbb{P} (M,\xi) \ar[d]^-{\mathbb{P} \tilde{\psi}}
	\\
		& \mathbb{P} (N,\xi_N).
	}\]
(Where an Engel local diffeomorphism $\mathbb{P} \tilde{\psi}$ is induced from $\tilde{\psi}$.)
\end{prop}
Let $(M_1, \xi_1)$ and $(M_2, \xi_2)$ be contact manifolds, and let $f : (M_1, \xi_1) \rightarrow (M_2, \xi_2)$ be a local contactomorphism.
Then, an Engel diffeomorphism $\mathbb{P}f : \mathbb{P}(M_1, \xi_1) \rightarrow \mathbb{P}(M_2, \xi_2)$ between the Cartan prolongations is induced.\par
Let ${\bf C}$ be a set of local contactomorphisms from $(M_1, \xi_1)$ to $(M_2, \xi_2)$.
Let ${\bf E}$ be a set of local Engel diffeomorphisms from $\mathbb{P}(M_1, \xi_1)$ to $\mathbb{P}(M_2, \xi_2)$.
The construction $(M, \xi) \mapsto \mathbb{P}(M, \xi)$ induces a map
\[
\mathbb{P} : {\bf C} \rightarrow {\bf E}
\]
\begin{cor}[\mbox{\cite{vogel2009existence}, cf. \cite[Corollary 1.20]{yamazaki2018engel}}] \label{corfullyfaithful}
The above map $\mathbb{P} : {\bf C} \rightarrow {\bf E}$ is a bijection.
\end{cor}
\begin{proof}
It is obvious because the counit of the adjunction $\mathbb{L} \dashv \mathbb{P}$ is an isomorphism.
(See Section \ref{sec1.4} and \cite{mac2013categories}.)
\end{proof}
The construction $(M, \xi) \mapsto \mathbb{P}(M, \xi)$ induces a group homomorphism $\operatorname{Aut}(M, \xi) \rightarrow \operatorname{Aut}(\mathbb{P}(M, \xi))$.
This homomorphism is an isomorphism for any contact manifold $(M, \xi)$ because of Corollary \ref{corfullyfaithful}.
(Corollary \ref{corfullyfaithful} means the functor $(M, \xi) \mapsto \mathbb{P}(M, \xi)$ is fully-faithful.
cf. Section \ref{sec1.4})
	\subsection{Engel manifolds and contact orbifolds} \label{sec1.3}
In this subsection, we generalize the discussion in the previous subsection to  a relationship between Engel manifolds and contact $3$-orbifolds.
First, we define a contact orbifold.
(See \cite{yamazaki2018engel} for a more detailed definition.)
\begin{defi}
Let $\Sigma$ be an orbifold.
A {\it contact structure} (resp. an {\it Engel structure}) on $\Sigma$ is a family of contact structures (resp. {\it Engel structures}) on each chart such that all local group actions and all transformation maps are contactomorphisms (resp. Engel diffeomorphisms).
The pair $(\Sigma, \xi)$ of the orbifold $\Sigma$ and a contact structure (resp. a Engel structure) $\xi$ on $\Sigma$ is called a {\it contact orbifold} (resp. an {\it Engel orbifold}). \par
Let $(\Sigma_1, \xi_1)$, $(\Sigma_2, \xi_2)$ be contact orbifolds (resp. Engel orbifolds).
A {\it local contactomorphism} (resp. a local Engel diffeomorphism) $f : (\Sigma_1, \xi_1) \rightarrow (\Sigma_2, \xi_2)$ is a local diffeomorphism $f : \Sigma_1 \rightarrow \Sigma_2$, which is a local contactomorphism (resp. a local Engel diffeomorphism) on each chart.
\end{defi}
As with Proposition \ref{proplspman}, the following proposition holds.
\begin{prop}[\mbox{cf. \cite[Proposition 2.4]{yamazaki2018engel}}] \label{proplsporb}
Let $(E, \mathcal{D})$ be an Engel manifold.
Suppose that the leaf space $\Sigma$ of the characteristic foiation of $(E, \mathcal{D})$ is an orbifold.
Then, the orbifold $\Sigma$ admits a contact structure.
\end{prop}
The pair of the leaf space of the characteristic foiation of the above Engel manifold $(E, \mathcal{D})$ and the above contact structure is a contact $3$-orbifold, which is denoted by $\mathbb{L}(E, \mathcal{D})$.
The Cartan prolongation of a contact $3$-orbifold is a patching of the Cartan prolongations of charts.
In general, the Cartan prolongation of a contact $3$-orbifold is an Engel orbifold.
\begin{defi} \label{goodchar}
Let $(E, \mathcal{D})$ be an Engel manifold.
Suppose that the leaf space of the characteristic foiation of $(E, \mathcal{D})$ is an orbifold.
(This assumption is not essential.
See Remark \ref{remnoneed}.)
We say that the Engel manifold $(E, \mathcal{D})$ has the {\it good characteristic foliation} if the Cartan prolongation $\mathbb{P}\mathbb{L}(E, \mathcal{D})$ of the contact orbifold $\mathbb{L}(E, \mathcal{D})$ is a manifold.
\end{defi}
Let $(E, \mathcal{D})$ be an Engel manifold with good characteristic foliation.
The developing map associated with $(E, \mathcal{D})$ is a patching of the developing maps associated with charts.
\begin{rem} \label{rembij2}
Let $\mathcal{L'}$ be the characteristic foliation of the Engel manifold $\mathbb{P}\mathbb{L}(E, \mathcal{D})$.
The leaf space of the Cartan prolongation $\mathbb{P}\mathbb{L}(E, \mathcal{D})$ is contactomorphic to the contact orbifold $\mathbb{L}(E, \mathcal{D})$.
Then, the developing map $\phi$ induces the following bijection $\phi_\ast$
(cf. Remark \ref{rembij1}).
\[
\{ \mbox{the leaves of the foliation $\mathcal{L}$} \} \overset{\phi_\ast}{\rightarrow} \{ \mbox{the leaves of the foliation $\mathcal{L'}$} \}.
\]
\end{rem}
An orbifold can be considered as a proper \'{e}tale Lie groupoid.
(See \cite{moerdijk2003introduction}.)
We will not discuss Lie groupoids in this paper, but we will give you some remarks.
In fact, the leaf space of a foliation is generally regarded as an (\'{e}tale) Lie groupoid by being identified with the holonomy groupoid.
The holonomy groupoid of the characteristic foliation of an Engel manifold has a contact structure.
(See \cite{yamazaki2018engel}.)
We write the holonomy groupoid with contact structure as $\mathbb{L}(E, \mathcal{D})$.
The Cartan prolongation of a $3$-dimendional contact \'{e}tale Lie groupoid is generally an Engel \'{e}tale Lie groupoid.
The necessary and sufficient condition that the Cartan prolongation of a contact \'{e}tale Lie groupoid is an Engel manifold is given by the following proposition.
\begin{prop}[\mbox{\cite[Corollary 4.5]{yamazaki2018engel}}]
The Cartan prolongation of a $3$-dimendional contact \'{e}tale Lie groupoid $\mathcal{G}$ is a manifold if and only if all of the following conditions are satisfied.
\begin{itemize}
\item $\mathcal{G}$ is an orbifold.
\item $\mathcal{G}$ is {\it positive}
(cf. \cite{yamazaki2018engel}).
\item $|\mathcal{G}_x|$ is odd for all $x \in \mathcal{G}_0$ where $\mathcal{G}_x$ is the isotoropy group at $x$.
\end{itemize}
\end{prop}
\begin{rem} \label{remnoneed}
If the Cartan prolongation $\mathbb{P}\mathcal{G}$ of a $3$-dimendional contact \'{e}tale Lie groupoid $\mathcal{G}$ is a manifold, then $\mathcal{G}$ is an orbifold.
That is why there is no need to discuss \'{e}tale Lie groupoids in this paper.
\end{rem}
\begin{rem} \label{remCartanReeb}
The orbifold $\mathcal{G}$ is the leaf space of the characteristic foliation of the Engel manifold $\mathbb{P}\mathcal{G}$.
In particular, the all holonomy groups of the characteristic foliation of $\mathbb{P}\mathcal{G}$ are finite.
The all leaves of the characteristic foliation of $\mathbb{P}\mathcal{G}$ are compact.
Then, the {\it Reeb local stability theorem} can be applied
(cf. Theorem \ref{reeb}).
\end{rem}
	\subsection{Categorical viewpoints} \label{sec1.4}
In this subsection, I will reconsider the previous discussion in terms of category theory.
Some concrete categories are defined as follows.
\begin{quote}
{\bf Contact} the category of contact $3$-manifolds. \\
${\bf Engel^t}$ the category of Engel manifolds with trivial characteristic foliation.
(See Definition \ref{trichar}.)\\
${\bf Contact'}$ the (higher) category of contact $3$-orbifolds whose Cartan prolongation is a manifold. \\
${\bf Engel^g}$ the category of Engel manifolds with good characteristic foliation.
(See Definition \ref{goodchar}.)
\end{quote}
The construction $(M, \xi) \mapsto \mathbb{P}(M, \xi)$ defines a functor $\mathbb{P} : {\bf Contact} \rightarrow {\bf Engel^t}$.
Let $\mathbb{L} : {\bf Engel^t} \rightarrow {\bf Contact}$ be a functor which sends each Engel manifold to the leaf space of the characteristic foliation.
Then, there exists an adjunction $\mathbb{L} \dashv \mathbb{P} : {\bf Engel^t} \rightarrow {\bf Contact^3}$.
The unit map $\phi : Id \rightarrow \mathbb{P} \circ \mathbb{L}$ of the adjunction $\mathbb{L} \dashv \mathbb{P}$ is the developing maps
(cf. Proposition \ref{propuniv}).
Corollary \ref{corfullyfaithful} means that the functor $\mathbb{P} : {\bf Contact} \rightarrow {\bf Engel^t}$ is fully faithful. \par
A developing map may not be an Engel diffeomorphism, but is a local Engel diffeomorphism.
Morphisms of the category ${\bf Engel^t}$ must be local Engel diffeomorphisms.
Then, morphisms of the category {\bf Contact} must be local contactomorphisms for the sake of the above argument.
(cf. Remark \ref{remlocal}) \par
Similarly, there exists an adjunction $\mathbb{L} \dashv \mathbb{P} : {\bf Engel^g} \rightarrow {\bf Contact'}$.
The unit map $\phi : Id \rightarrow \mathbb{P} \circ \mathbb{L}$ of the adjunction $\mathbb{L} \dashv \mathbb{P}$ is the developing map. \par
Remark that the class ${\bf Contact'}$ may not be a category, but is a higher category.
(If an orbifold is defined as a Lie groupoid, then ${\bf Contact'}$ is not a category, but is a $2$-category.)
Then, $\mathbb{L} : {\bf Engel^g} \rightarrow {\bf Contact'}$ and $\mathbb{P} : {\bf Contact'} \rightarrow {\bf Engel^g}$ may not be functors, but are lax functors.
However, $\mathbb{P} \circ \mathbb{L} : {\bf Engel^g} \rightarrow {\bf Engel^g}$ is a functor, because it is a lax functor between ($1$-)categories.
Therefore, we only need to discuss ($1$-)functors in this paper.
	\section{The main result} \label{sec2}
We will show the following theorem.
\begin{thm} \label{main}
Let $(E, \mathcal{D})$ be a connected Engel manifold with good characteristic foliation
(cf. Definition \ref{goodchar}).
If the developing map $\phi$ associated with $(E, \mathcal{D})$ is not a covering map, then a group homomorphism $\Phi = \mathbb{P} \circ \mathbb{L} : \operatorname{Aut}(E, \mathcal{D}) \rightarrow \operatorname{Aut}(\mathbb{P}\mathbb{L}(E, \mathcal{D}))$ is injective.
\end{thm}
	\subsection{Foliations and lemmas}
We will use the following theorem without proof
(cf. \cite{tamura2006topology}).
\begin{thm}[Reeb local stability theorem] \label{reeb}
Let $P$ be a manifold, and let $\mathcal{F}$ be a foliation on $P$.
Suppose that a leaf $L$ of $\mathcal{F}$ is compact, and that the holonomy group at a point $x \in L$ is finite.
Then, for any open neighborhood $U' \subset P$ of $L$, there exists a tubular neighborhood $U \subset P$ of $L$ with a projection $p : U \rightarrow L$ such that
\begin{itemize}
\item $U \subset U'$,
\item the open set $U$ is a union of some compact leaves of $\mathcal{F}$,
\item a fiber $p^{-1} (z)$ intersects transversally with $\mathcal{F}$ for any point $z \in L$.
\end{itemize}
\end{thm}
\begin{rem} \label{tubesec}
Fix any Riemannian metric on the manifold $P$.
A tubular neighborhood of the leaf $L$ is constructed by starting geodesics from $L$ to normal direction. \par
Suppose that there exists a slice $S \subset P$ of $\mathcal{F}$.
(i.e. The submanifold $S \subset P$ intersects transversally with the foliation $\mathcal{F}$, and the dimension of the manifold $S$ is $\operatorname{dim}(P) - \operatorname{dim}(\mathcal{F})$.)
Take a point $x \in S \cap L$.
There is a locally defined metric $g$ around $x \in P$ such that the slice $S$ is totally geodesic around $x \in S$ in local.
The domain of the metric $g$ can be extended to the whole domain $P$ by a partition of unity.
Then, the tubular neighborhood $U$ of the leaf $L$ with a projection $p : U \rightarrow L$ satisfies $p^{-1}(x) \subset S$.
This discussion will be used later.
\end{rem}
Moreover, we will prepare some technical lemmas for the proof of Theorem \ref{main}.
\begin{lem} \label{reeb2}
Suppose that $P$, $\mathcal{F}$, and $L$ and the holonomy is the same as in the {\it Reeb local stability theorem}
(Theorem \ref{reeb}).
Moreover, suppose that the dimension of the foliation $\mathcal{F}$ is $1$.
Let $U \subset P$ be any tubular neighborhood of $L$ with a projection $p : U \rightarrow L$ as in the {\it Reeb local stability theorem}.
Then, there exists a complete non-singular vector field $X$ on $U$ tangent to $\mathcal{F}$ such that $\rho_t (p^{-1} (z)) \subset p^{-1} (z)$ if and only if $t \in \mathbb{Z}$ for any point $z \in L$, where $\rho_t = \operatorname{Exp}(tX) : U \rightarrow U$ is the flow of $X$.
\end{lem}
\begin{proof}
Because $L \cong S^1$, there exists a complete non-singular vector field $Y$ on the leaf $L$ such that $\operatorname{Exp}(tY) = id$ if and only if $t \in \mathbb{Z}$.
We regard the foliation $\mathcal{F}$ as an integrable distribution.
Then, $\mathcal{F}$ is a vector bundle $\mathcal{F} \rightarrow P$ on $P$.
The following diagram is a fiber product diagram because a fiber $p^{-1} (z)$ and the foliation $\mathcal{F}$ intersect transversally for any $z \in L$:
	\[\xymatrix{
		\mathcal{F} \ar[r]^-{dp} \ar[d]
		& TL \ar[d]
	\\
		U \ar[r]_-{p}
		& L.
	}\]
Let $X$ be the pullback of $Y$ along the map $p$.
The vector field $X$ is tangent to $\mathcal{F}$.
The open set $U$ is a union of some compact leaves of $\mathcal{F}$.
The vector field $X$ is complete because all leaves in the open set $U$ are compact.
The vector field $X$ is non-singular because $Y$ is non-singular.
Take any point $z \in L$.
For any point $w \in p^{-1} (z)$,
\[
\begin{array}{lcl}
\rho_t (w) \in p^{-1} (z) & \Leftrightarrow & \rho_t (p(w)) = z \\
& \Leftrightarrow & t \in \mathbb{Z}. 
\end{array}
\]
Therefore, $\rho_t (p^{-1} (z)) \subset p^{-1} (z)$ if and only if $t \in \mathbb{Z}$.
\end{proof}
\begin{lem} \label{noncov}
Let $\psi : \mathbb{R} \rightarrow S^1$ be a local diffeomorphism. 
There exists an open embedding $\tau : \mathbb{R} \rightarrow \mathbb{R}$ such that the following diagram is commutative, where $\pi(\theta) = [\theta]$:
	\[\xymatrix{
		& \mathbb{R} \ar@{->>}[d]^-{\pi}
	\\
		\mathbb{R} \ar[ur]^-{\tau} \ar[r]_-{\psi}
		& \mathbb{R} /2 \pi \mathbb{Z} \cong S^1.
	}\]
\end{lem}
\begin{proof}
The affine space $\mathbb{R}$ is simply connected.
There exists a lift $\tau : \mathbb{R} \rightarrow \mathbb{R}$ of the local diffeomorphism $\psi$ for the universal covering $\pi : \mathbb{R} \twoheadrightarrow \mathbb{R} /2 \pi \mathbb{Z} \cong S^1$.
Then, the local diffeomorphism $\tau$ is monotonic.
The map $\tau$ is an open embedding.
\end{proof}
	\subsection{Proof of the main result}
First, we prepare several symbols.
\begin{defi} \label{closure}
Let $(E, \mathcal{D})$ be an Engel manifold, and let $\mathcal{L}$ be the characteristic foliation of $(E, \mathcal{D})$.
For any point $e \in E$, denote a leaf of the foliation $\mathcal{L}$ including $e \in E$ as $\widehat{e}$.
For any subset $A \subset E$, define a subset $\widehat{A} \subset E$ as $\widehat{A} = \bigcup_{e \in A} \widehat{e}$.
\end{defi}
Let $(E, \mathcal{D})$ be a connected Engel manifold with good characteristic foliation $\mathcal{L}$, and let $\phi : E \rightarrow \mathbb{P}\mathbb{L}(E, \mathcal{D})$ be the developing map of $(E, \mathcal{D})$.
Theorem \ref{main} follows immediately from the following two lemmas.
\begin{lem} \label{lem1}
If there exists a leaf $L$ of the characteristic foliation $\mathcal{L}$ such that the map $\phi |_L : L \rightarrow \widehat{\phi(L)}$ (cf. Definition \ref{closure}) is not a covering map, then a group morphism $\Phi = \mathbb{P} \circ \mathbb{L} : \operatorname{Aut}(E, \mathcal{D}) \rightarrow \operatorname{Aut}(\mathbb{P}\mathbb{L}(E, \mathcal{D}))$ is injective.
\end{lem}
\begin{lem} \label{lem2}
If the map $\phi |_L : L \rightarrow \widehat{\phi(L)}$ (cf. Definition \ref{closure}) is a covering map for any leaf $L$ of the characteristic foliation $\mathcal{L}$, then the developing map $\phi : E \rightarrow \mathbb{P}\mathbb{L}(E, \mathcal{D})$ is a covering map.
\end{lem}
\begin{proof}[Proof of Theorem \ref{main}]
Suppose that the developing map $\phi$ is not a covering map.
By the counterpart of Lemma \ref{lem2}, there exists a leaf $L$ of the characteristic foliation $\mathcal{L}$ such that the map $\phi |_L : L \rightarrow \widehat{\phi(L)}$ is not a covering map.
Then, the group homomorphism $\Phi$ is injective by Lemma \ref{lem1}.
\end{proof}
We prove these lemmas in turn.
\begin{proof}[Proof of Lemma \ref{lem1}]
Suppose $\Phi (f) = \Phi (g)$ for two Engel diffeomorphisms $f, g \in \operatorname{Aut}(E, \mathcal{D})$. 
Let $\sigma$ be the local Engel diffeomorphism $\Phi (f) (= \Phi (g))$, and let $L' = \widehat{\phi(L)}$. 
We discuss the following diagram:
	\[\xymatrix{
		\phi^{-1} (\operatorname{Im}(\sigma |_{L'})) \ar[r]^-{\phi}
		& \operatorname{Im}(\sigma |_{L'})
	\\
		L \ar@<0.5ex>[u]^-{f} \ar@<-0.5ex>[u]_-{g} \ar[r]_-{\phi}
		& L' \ar[u]^-{\sigma}.
	}\]
$\phi^{-1} (\operatorname{Im}(\sigma |_{L'}))$ is a leaf of the foliation $\mathcal{L}$ because of Remark \ref{rembij2}.
Subsets $f(L)$ and $g(L)$ are leaves of the foliation $\mathcal{L}$ because the maps $f$ and $g$ are Engel diffeomorphisms.
We obtain an equality $\phi^{-1} (\operatorname{Im}(\sigma |_{L'})) =f(L) = g(L)$.
The maps $\phi |_L$ and $\phi |_{\phi^{-1} (\operatorname{Im}(\sigma |_{L'}))}$ are not covering maps.
Both maps $\phi |_L$ and $\phi |_{\phi^{-1} (\operatorname{Im}(\sigma |_{L'}))}$ can be regarded as local diffeomorphisms $\mathbb{R} \rightarrow S^1$. 
By Lemma \ref{noncov}, the maps $\phi |_L$ and $\phi |_{\phi^{-1} (\operatorname{Im}(\sigma |_{L'}))}$ may be regarded as the map $(a, b), (c, d) \subsetneq \mathbb{R} \twoheadrightarrow \mathbb{R} /2 \pi \mathbb{Z} \cong S^1;\theta \mapsto [\theta]$ respectively.
Reversing the orientation, if necessary, we can assume $b, d \neq \infty$. 
We discuss the following diagram:
	\[\xymatrix{
		(c, d) \ar[r]^-{\phi}
		& S^1
	\\
		(a, b) \ar@<0.5ex>[u]^-f \ar@<-0.5ex>[u]_-g \ar[r]_-{\phi}
		& S^1 \ar[u]^-{\sigma}.
	}\]
A congruence equality $f(t) \equiv g(t)$ mod $2 \pi$ holds for any point $t \in (a, b)$.
A function $f |_L -g |_L$ is a constant function because of a differential equality $f' = g'$.
Then, there exists a integer $k$ such that a equality $f |_L - g |_L = 2k \pi$ holds.
Take a sequence $\{ t_n \} \subset (a, b)$ converging to the point $b$.
An equality $f(t_n) = g(t_n) + 2k \pi$ holds.
Both sequences $\{ f(t_n) \}$ and $\{ g(t_n) \}$ converge to the point $d$ because a sequence $\{ \sigma([t_n]) (= \phi(f(t_n)) = \phi(g(t_n))) \}$ converge to the point $[d] \in S^1$ and the maps $f : (a, b) \rightarrow (c, d)$ and $g : (a, b) \rightarrow (c, d)$ are diffeomorphisms.
An equality $d = d + 2k \pi$ holds.
We obtain a solution $k = 0$.
Then, we obtain an equality $f |_L = g |_L$.
The two lifts $f$ and $g$ of a map $\sigma \circ \phi$ are identical on the connected domain $E$. 
Therefore, the group morphism $\Phi$ is injective.
\end{proof}
\begin{proof}[Proof of Lemma \ref{lem2}]
Denote the characteristic foliation of the Engel manifold $\mathbb{P}\mathbb{L}(E, \mathcal{D})$ as $\mathcal{L'}$. \par
Suppose that the map $\phi |_L : L \rightarrow \widehat{\phi(L)}$ is a covering map for any leaf $L$ of the characteristic foliation $\mathcal{L}$.
We will show that the developing map $\phi$ is a covering map. \par
Let $C \subset E$ and $N \subset E$ be the sets as the following.
\[
\begin{array}{lcl}
C & = & \{ e \in E \, | \, \mbox{The leaf $\widehat{e}$ of the foliation $\mathcal{L}$ is compact.} \}, \\
N & = & \{ e \in E \, | \, \mbox{The leaf $\widehat{e}$ of the foliation $\mathcal{L}$ is not compact.} \}.
\end{array}
\]
Remark $E = C \coprod N$.
The subset $C \subset E$ is open because of the {\it Reeb local stability theorem} (cf. Theorem \ref{reeb}).
(See Remark \ref{remCartanReeb} for reasons why the {\it Reeb local stability theorem} can be applied.) \par
Let $\{ C_i \}_i$ be the set of connected components of $C$.
Then, the map $\phi |_{C_i} : C_i \rightarrow \widehat{\phi(C_i)}$ is proper.
Therefore, the map $\phi |_{C_i}$ is a finite covering map. \par
We will show that the subset $C_i \subset E$ is closed and open in the whole space $E$.
It is obvious that the subset $C_i \subset E$ is open in the whole space $E$. 
We have only to show that the subset $C_i \subset E$ is closed in the whole space $E$.
Take any point $x$ in the closure of the subset $C_i \subset E$ in the whole space $E$, and let $y = \phi (x)$.
We will show that the point $x$ is included in the subspace $C_i$.
We have only to show that the point $x$ is included in the subspace $C$ because the subset $C_i \subset C$ is closed in the subspace $C$.
Let $L = \widehat{x}$, and let $L' = \widehat{y}$.
(Remark, the leaf $L$ of the foliation $\mathcal{L}$ may be not compact, but the leaf $L'$ of the foliation $\mathcal{L'}$ is compact.)
Take a slice $O \subset E$ of $\mathcal{L}$ including the point $x$ such that the map $\phi |_O$ is an embedding,.
(There exists such a slice $O$ because the map $\phi$ is a local diffeomorphism.)
By Lemma \ref{reeb2}, there exist
\begin{itemize}
\item a tubular neighborhood $V \subset \mathbb{P}\mathbb{L}(E, \mathcal{D})$ of the leaf $L'$ with a projection $q : V \rightarrow L'$ and
\item a complete non-singular vector field $Y$ on the open set $V$ tangent to the foliation $\mathcal{L'}$
\end{itemize}
such that these satisfy the conditions in Lemma \ref{reeb2}.
Specifically these satisfy the following conditions.
\begin{itemize}
\item $V = \widehat{V}$,
\item a fiber $q^{-1} (z)$ intersects transversally with the foliation $\mathcal{L'}$ for any point $z \in L'$,
\item $\rho_t (q^{-1} (z)) \subset q^{-1} (z)$ if and only if $t \in \mathbb{Z}$ for any point $z \in L'$, where $\rho_t = \operatorname{Exp}(tY) : V \rightarrow V$ is the flow of the vector field $Y$.
\end{itemize}
Furthermore, by Remark \ref{tubesec}, we may assume the following.
\begin{itemize}
\item $q^{-1}(y) \subset \phi(O)$.
\end{itemize}
Let $S = (\phi |_O)^{-1}(q^{-1}(y)) (\subset O)$.
(Then, the map $\phi$ induces a diffeomorphism $S \cong q^{-1}(y)$.)
Let $U = \phi^{-1}(V)$, and let $X$ be a pullback of the vector field $Y$ by the developing map $\phi$.
(The vector field $X$ is a non-singular vector field on the open set $U \subset E$.)
The vector field $X$ is complete because the map $\phi |_{L_0}$ is a covering map for any leaf $L_0$ of the foliation $\mathcal{L}$.
Let $\tilde{\rho}_t = \operatorname{Exp}(t X)$ be the flow of the vector field $X$.
Then, a map $\tilde{\rho} : S \times \mathbb{R} \rightarrow \tilde{U} ; (w,t) \mapsto \tilde{\rho}_t (w)$ is a surjective local diffeomorphism.
There exists a sequence $\{ x_n \} \subset C_i$ converging to the point $x$, by definition of the point $x$.
Moreover, we can assume that there exists a sequence $\{ (w_n, t_n) \} \subset S \times \mathbb{R}$ such that equalitys $x_n = \tilde{\rho}_{t_n}(w_n)$, because the map $\tilde{\rho}$ is a local diffeomorphism.
The sequence $\{ (w_n, t_n) \}$ converges to the point $(x, 0)$ because of an equality $x = \tilde{\rho}_0(x)$.
Then, the sequence $\{ w_n \}$ converges to the point $x$.
The point $w_n$ is included in the set $C_i$ because of a inclusion relation $\widehat{w_n} = \widehat{x_n} \subset C_i$.
Let $k$ be the order of the holonomy group of the leaf $\widehat{y}$ of the foliation $\mathcal{L'}$. 
Then, we obtain an equality $w_n = \tilde{\rho}_{k n_i} (w_n)$, where the integer $n_i$ is the degree of the finite covering map $\phi |_{C_i} : C_i \rightarrow \widehat{\phi(C_i)}$.
If we take the limit with respect to the index $n$ on both sides, we get an equality $x = \tilde{\rho}_{k n_i} (x)$. 
This means that a trajectory of the vector field $X$ through the point $x$ is periodic.
In other words, the subset $\widehat{x} \subset E$ is compact.
Therefore, the point $x$ is included in the subspace $C$.
The component $C_i$ of the subspace $C$ is closed in the subspace $C$.
The point $x$ is included in the component $C_i$.
Therefore, the subset $C_i \subset E$ is closed in the whole space $E$. \par
We will show that the developing map $\phi$ is a covering map.
\begin{description}
  \item[Case.1] $C \neq \emptyset$.
\end{description}\par
Then, we obtain an equality $E=C=C_i$ ($^{\exists!} i$) because the manifold $E$ is connected and the subset $C_i \subset E$ is closed and open. 
The developing map $\phi$ is proper.
The developing map $\phi$ is a finite covering map.
\begin{description}
  \item[Case.2] $C = \emptyset$.
\end{description}\par
Then, we obtain an equality $E=N$. 
In particular, all leaves of the foliation $\mathcal{L}$ are simply connected. 
All of holonomy groups are trivial. 
(Then, the leaf space $\mathbb{L}(E, \mathcal{D})$ is a manifold.)
The developing map $\phi$ is surjective because the map $\phi |_L : L \rightarrow \widehat{\phi(L)}$ is surjective for any leaf $L$ of the foliation $\mathcal{L}$
(cf. Remark \ref{rembij1}).
Take any point $y$ in the Cartan prolongation $\mathbb{P}\mathbb{L}(E, \mathcal{D})$, and fix a point $x \in \phi^{-1} (y)$.
Let $L = \widehat{x}$, and let $L' = \widehat{y}$.
Same to the above, take $V$, $q$, $Y$, $\rho$, $U$, $X$, $\tilde{\rho}$, $S$.
The following diagram is commutative:
	\[\xymatrix{
		S \times \mathbb{R} \ar[d]_-{\tilde{\rho}}^-{\cong} \ar[r]^-{\phi \times id}_-{\cong}
		& q^{-1} (y) \times \mathbb{R} \ar[d]^-{\rho}
	\\
		U \ar[r]^-{\phi}
		& V.
	}\]
We will explain that the map $\tilde{\rho}$ is a diffeomorphism.
The map $\tilde{\rho}$ is obviously a local diffeomorphism by its construction.
We have to show that the map $\tilde{\rho}$ is bijective.
We will construct an inverse map of the map $\tilde{\rho}$.
Take any point $z \in U$.
There exists a unique point $w \in S$ such that $\widehat{\phi(z)} \cap S = \{ w \}$, because the holonomy group at the point $\phi(z)$ is trivial.
Then, an equality $\widehat{\phi(z)} = \widehat{\phi(w)}$ holds.
We obtain an equality $\widehat{z} = \widehat{w}$ because of the correspondence $\phi_\ast$ in Remark \ref{rembij1}.
There exists a point $t \in \mathbb{R}$ such that $\tilde{\rho}_t (w) = z$.
Moreover, the point $t$ is unique because the vecor field $X$ is not periodic. 
(If it is periodic, then the foliation $\mathcal{L}$ has a compact leaf.
This is in contradiction to $C = \emptyset$.)
We define a map by $z \mapsto (w, t)$.
The map $z \mapsto (w, t)$ is the inverse map of the map $\tilde{\rho}$.
The map $\tilde{\rho}$ is bijective.
Therefore, the map $\tilde{\rho}$ is a diffeomorphism. \par
The map $\rho : p^{-1} (y) \times \mathbb{R} \rightarrow V ; (z,t) \mapsto \rho_t (z)$ is obviously a covering map.
Then, the developing map $\phi$ is a covering map in local around the point $y$.
Therefore, the developing map $\phi$ is a covering map in global, because the point $y$ was taken arbitarily.
\end{proof}
	\subsection{Application: An answer of AIM problem}
Let $(E, \mathcal{D})$ be a connected Engel manifold with trivial characteristic foliation $\mathcal{L}$
(cf. Definition \ref{trichar}).
Let $\phi$ be the developing map. \par
Define the function $\sigma : \mathbb{L}(E, \mathcal{D}) \rightarrow \mathbb{Z}_{\ge 0} \cup \{ \infty \}$ as $\sigma(L) = \min \{ \# \phi^{-1} (y) \, | \, y \in \widehat{\phi(L)} \}$. 
(For a set $A$, a number $\# A$ is a cardinal of the set $A$.
If the set $A$ is a infinite set, then we define $\# A = \infty$.)
We call the function $\sigma$ the {\it twisting number function}.
\begin{prop} \label{prop}
For any Engel diffeomorphism $f \in \operatorname{Aut}(E,\mathcal{D})$, the induced diffeomorphism $\mathbb{L}f : \mathbb{L}(E, \mathcal{D}) \rightarrow \mathbb{L}(E, \mathcal{D})$ preserves the twisting number function $\sigma$.
\end{prop}
\begin{proof}
Let $L$ be a point of the leaf space $\mathbb{L}(E, \mathcal{D})$.
$L$ is a leaf of the foliation $\mathcal{L}$. 
We consider the following diagram:
	\[\xymatrix{
		f(L) \ar[r]^-{\phi}
		& \widehat{\phi(f(L))}
	\\
		L \ar[u]^-{f}_-{\cong} \ar[r]^-{\phi}
		& \widehat{\phi(L)} \ar[u]^-{\mathbb{P}\mathbb{L}f}_-{\cong}.
	}\]
The maps $\phi |_L : L \rightarrow \widehat{\phi(L)}$ and $\phi |_{f(L)} : f(L) \rightarrow \widehat{\phi(f(L))}$ are isomorphic.
Then, we obtain an equality $\min \{ \# \phi^{-1} (y) \, | \, y \in \widehat{\phi(L)} \} = \min \{ \# \phi^{-1} (y) \, | \, y \in \widehat{\phi(f(L))} \}$.
This means that an equality $\sigma(L) = \sigma(\mathbb{L}f(L))$ holds.
\end{proof}
\vspace{10pt}
The following example is the universal covering of the Cartan prolongation of the affine space $\mathbb{R}^3$ with standard contact structure.
Then, the twisting number function is trivial.
\begin{ex} \label{ex1}
Let $E = \mathbb{R}^4$.
Denote its coordinates as $(x,y,z,\theta) \in \mathbb{R}^4$.
Let $\mathcal{D} = \langle \partial_\theta , \cos(\theta /2) X + \sin(\theta /2) Y \rangle$, (for $X = \partial_x - y \partial_z$, $Y = \partial_y$). 
The pair $(E, \mathcal{D})$ is an Engel manifold. \par
The distribution $\mathcal{L} (= \langle \partial_\theta \rangle)$ is the characteristic foliation.
The distribution $\mathcal{E} (= \mathcal{D}^2)$ is represented as $\mathcal{E} = \langle \partial_\theta, X, Y \rangle$. 
The leaf space $M$ of the foliation $\mathcal{L}$ is represented as $M \cong \mathbb{R}^3 \ni (x,y,z)$, and the induced contact structure $\xi (= \mathcal{E}/\mathcal{L})$ on $M$ is represented as $\xi \cong \langle X, Y \rangle$. 
We have an identification $M \times S^1 \cong \mathbb{P}(\xi) ; (\bm{x}, [\theta]) \mapsto (\langle \cos(\theta /2) X_{\bm{x}} + \sin(\theta /2) Y_{\bm{x}} \rangle, \bm{x})$. 
A map $\phi : E \rightarrow M \times S^1 ; (\bm{x}, \theta) \mapsto (\bm{x}, [\theta])$ is identified with the developing map.
Then, the twisting number function $\sigma$ is a constant function $\sigma \equiv \infty$.
\end{ex}
\vspace{10pt}
We slightly modify Example \ref{ex1}.
\begin{ex}
Fix a point $\bm{x}_0 \in \mathbb{R}^3$ and an integer $n \in \mathbb{Z}_{\ge 0}$. \par
Let $E = \mathbb{R}^4 - \{ \bm{x}_0 \} \times ((-\infty,-n\pi] \cup [n\pi +\epsilon,\infty))$, for some constant $\epsilon \in (0,2\pi]$.
The subset $E \subset \mathbb{R}^4$ is open.
Then, the subspace $E \subset \mathbb{R}^4$ is an open submanifold.
Let $\mathcal{D}$ be as in Example \ref{ex1}.
The pair $(E, \mathcal{D})$ is an Engel manifold. \par
A map $\phi : E \rightarrow \mathbb{R}^3 \times S^1 ; (\bm{x}, \theta) \mapsto (\bm{x}, [\theta])$ is identified with the developing map.
Then, the twisting number function $\sigma : \mathbb{R}^3 \rightarrow \mathbb{Z}_{\ge 0} \cup \{ \infty \}$ is the following.
\[
  \sigma(\bm{x}) = \begin{cases}
    n & (\bm{x} = \bm{x}_0) \\
    \infty & (\mbox{otherwise}).
  \end{cases}
\]
\end{ex}
\vspace{10pt}
We further modify these examples to make an Engel manifold whose automorphism group is trivial.
This implies that we give a positive answer to the AIM Problem.
\begin{ex} \label{ex3}
Take a countable dense subset $Q = \{ \bm{x}_n \}_{n=1}^\infty \subset \mathbb{R}^3$. 
Let $E = \mathbb{R}^4 - \bigcup_{n = 1}^{\infty} \{ \bm{x}_n \} \times ((-\infty,-n\pi] \cup [n\pi +\epsilon,\infty))$, for some constant $\epsilon \in (0,2\pi])$.
The subset $E \subset \mathbb{R}^4$ is open.
Then, the subspace $E \subset \mathbb{R}^4$ is an open submanifold.
Let $\mathcal{D}$ be as in Example \ref{ex1}.
The pair $(E, \mathcal{D})$ is an Engel manifold. \par
A map $\phi : E \rightarrow \mathbb{R}^3 \times S^1 ; (\bm{x}, \theta) \mapsto (\bm{x}, [\theta])$ is identified with the developing map.
Then, the twisting number function $\sigma : \mathbb{R}^3 \rightarrow \mathbb{Z}_{\ge 0} \cup \{ \infty \}$ is the following.
\[
  \sigma(\bm{x}) = \begin{cases}
    n & (\bm{x} = \bm{x}_n) \\
    \infty & (\mbox{otherwise}).
  \end{cases}
\]
\end{ex}
The author had come up with this construction by himself, but he mistakenly thought that the above $E$ would not be a manifold.
When the author introduced this AIM problem and our main result to Mitsumatsu, he proposed the same construction independently of the author, which led to the final solution.
\begin{thm}[Mitsumatsu, Y.] \label{app}
Let $(E, \mathcal{D})$ be the Engel manifold as in Example \ref{ex3}.
The autorphism group $\operatorname{Aut}(E, \mathcal{D})$ is trivial.
\end{thm}
\begin{proof}
For any Engel diffeomorphism $f \in \operatorname{Aut}(E,\mathcal{D})$, the induced diffeomorphism $\mathbb{L}f : \mathbb{R}^3 \rightarrow \mathbb{R}^3$ is the identity on $Q$ by Proposition \ref{prop}.
The diffeomorphism $\mathbb{L}f$ is the identity on the leaf space $\mathbb{R}^3$ because the subset $Q \subset \mathbb{R}^3$ is dense.
The Engel diffeomorphism $f$ is the identity by Theorem \ref{main} (or Lemma \ref{lem1}).
\end{proof}

\bibliography{engel,contact,symplectic,category,homotopy,foliation,Lie_groupoid}
\bibliographystyle{plain}

Koji Yamazaki\\
Department of Mathematics, Tokyo Institute of Technology\\
2-12-1, Ookayama Meguro-ku, Tokyo, Japan\\
e-mail:yamazaki.k.as@m.titech.ac.jp

\end{document}